\documentclass[11pt]{article}
\usepackage[margin=1.08in]{geometry}
\usepackage{fontspec}
\newfontfamily\greekfont[Script=Greek]{FreeSerif.otf}
\usepackage{amsmath,amssymb}
\usepackage{microtype}
\usepackage{enumitem}
\usepackage{authblk}
\usepackage{csquotes}
\usepackage{xcolor}
\usepackage{hyperref}
\hypersetup{colorlinks=true,linkcolor=blue!45!black,citecolor=blue!45!black,urlcolor=blue!45!black}
\usepackage[backend=biber,style=numeric,sorting=nyt,doi=true,url=true,isbn=false]{biblatex}
\newcommand{\gnothi}{{\greekfont γνῶθι σεαυτόν}}
\newcommand{\practice}{\mathcal P}
\newcommand{\concept}{\mathcal C}
\newcommand{\enr}{\mathrm{Enr}}
\newcommand{\moves}{\mathcal{A}}
\newcommand{\negative}{\mathcal{N}}
\title{The Work of Mathematical Production in the Age of AI}
\author[1]{Asaf Hadari}
\author[2]{Thomas Koberda}

\affil[1]{Department of Mathematics,
University of Hawai`i at Manoa,
Honolulu, HI 96822, USA}

\affil[2]{Department of Mathematics,
University of Virginia,
Charlottesville, VA 22904, USA}

\affil[1]{\texttt{hadari@hawaii.edu}}
\affil[2]{\texttt{thomas.koberda@virginia.edu}}

\date{\today}
\begin{document}
\maketitle

\begin{abstract}
The emergence of powerful artificial intelligence has unsettled the relationship between mathematical production and the human intellectual activity that mathematical artifacts have traditionally documented. We argue that mathematical practice cannot be adequately understood through the production of mathematical knowledge alone, and distinguish four interrelated dimensions of inquiry: knowledge, understanding, making sense, and figuring out.

Building on Kitcher, Dewey, and subsequent philosophy of mathematical practice, we develop an expanded framework incorporating the cognitive and computational capacities of mathematical communities, their conceptual organizations, provisional investigative models, and the role of unsuccessful inquiry. We propose reconceptualizing artificial intelligence as a manifestation of physically organized computation whose mathematical behavior itself can become an object of experimental investigation.

This perspective motivates a program of \emph{natural mathematics}, by which we mean the construction and comparison of artificial intellectual communities developing mathematical practices under different cognitive, conceptual, and normative conditions. Our aim is not simply to automate mathematical research or increase the production of theorems, but to investigate the robustness and contingency of mathematical concepts, the possibilities of alternative mathematical practices, and what such investigations might reveal about human mathematical cognition.

Finally, we consider the future of human intellectual agency in a world where mathematical production need no longer depend upon humans. We argue for intergenerational epistemic stewardship as a justification for sustaining human mathematical traditions and opportunities for intellectual participation, independently of humanity's continued indispensability to mathematical progress.
\end{abstract}

\section{Introduction}

It is by now obvious to practitioners of mathematics that the ``aura'' of mathematical papers, the primary artifacts of mathematical intellectual activity, has come into question. Walter Benjamin \cite{Benjamin2008} famously connected the aura of a work of art with its singular existence, authenticity, historical situation, and relationship to tradition. The aura of a mathematical paper is not identical with that of a work of art; one might reasonably define it not through uniqueness, but through the paper's relationship to an identifiable act of mathematical inquiry. The diverse posts on \emph{Proofs and Prompts} \cite{ProofsPrompts} suggest that this relationship has suddenly become unstable, and many mathematicians are desperately trying to preserve the aura.

Every new preprint either consciously makes use of generative artificial intelligence or inhabits a world in which artificial intelligence is a prominent feature. Even a paper that is honestly ``AI vegan'' raises questions that its authors did not intend, that inevitably diminish the romantic aura of the paper, and that until recently would not have been posed: Is it possible to be fully free of the influence of artificial intelligence, and what bearing does the answer have on the paper itself? Could artificial intelligence have done this, and perhaps done it better? If certain tasks attendant upon mathematical production are indisputably done better by AI, and even augment human-centered goals such as the human understanding of mathematics, is it a good use of resources, financial or intellectual, to write AI-vegan papers?

In this essay, we investigate the relationship between the production of mathematical knowledge and the human activity of doing mathematics, and the changes in that relationship brought into sharp relief by powerful artificial intelligences. We will confront four questions, in an order that takes us from mathematical practice to new possibilities for inquiry, and finally to the place of humans within it:
\begin{enumerate}[leftmargin=2.2em]
\item What is the structure of human mathematical practice, and how precisely has artificial intelligence affected it?
\item How should mathematicians conceive of artificial intelligence, and how should this change mathematical inquiry?
\item What opportunities does artificial intelligence make available beyond the generation of more and more theorems---opportunities hitherto inaccessible?
\item What is the role of the human mathematician, and more generally the human intellectual agent, in a world with powerful artificial intelligence?
\end{enumerate}

Our argument begins by distinguishing mathematical
knowledge from the further activities of understanding,
making sense, and figuring out. We then propose a
naturalistic conception of artificial intelligence
as physically organized computation whose mathematical
behavior can be investigated experimentally.
This motivates our central proposal: to construct
and compare artificial intellectual communities
developing alternative mathematical practices,
and thereby investigate the contingency and
robustness of mathematical conceptual organization.
Our aim is to broaden mathematical inquiry and what it might reveal about human thought.
Finally, we consider the implications for human
intellectual agency, and argue that the possibility
of mathematics without human participation does
not eliminate the value of sustaining human
intellectual traditions, particularly as a form
of intergenerational epistemic stewardship.

\section{The Uncanny Mathematical Paper}

In late September 2026, one of us was asked to provide a brief opinion for a journal on a paper by a prominent author. As he read, a suspicion grew that the paper had been entirely generated by AI. He did not feel deceived: indeed, at the bottom there was a disclosure that the paper had been prepared with the help of AI, though this did not rule out its having been entirely generated by AI. An obvious possibility was that the author had an idea, prompted an AI to develop it and write it up, and then read the paper line by line to certify its correctness. The disclosure was too vague to permit a more precise conjecture.

It was difficult for him to understand exactly what he found unsettling or uncanny about the paper. He is not opposed to using artificial intelligence in mathematical research, and has spent considerable time exploring its possibilities. The usual objections---plagiarism, inauthenticity, fear of artificial intelligence itself, or doubts about the correctness of its results---did not apply, or at least did not explain his unease, nor why he instinctively felt that this apparently AI-generated paper ought to be held to different standards from an AI-vegan paper produced a year earlier.

He had long thought of a mathematical paper as a declaration of the notability of a mathematical investigation and of its worthiness of others' attention, not simply of its material conclusions. The mathematical artifact traditionally served as documentation of an epistemic act. The relationship between an author's epistemic acts and the mathematical artifact is irrevocably altered once artificial intelligence becomes part of the environment in which the author works. It therefore seemed natural that the standards by which artifacts of mathematical production are evaluated could be determined only after clarifying the function of the artifact itself.

Further reflection suggested that this perspective was too coarse. The traditional understanding that mathematical artifacts were proxies for epistemic acts was more of a sociological convention than an essential relationship. A fully human-authored paper, painstakingly developed over many years with deep understanding, might communicate the resulting conceptual organization only partially. Another human-authored paper might result from mechanical manipulations carefully checked and documented with little attention paid to significance or conceptual understanding. A third, AI-generated result might be deeply studied and understood by a human, its significance carefully and thoroughly explained. All three possibilities describe familiar occurrences in contemporary mathematical production, and in each case the relationship between artifact and epistemic activity differs. It would be too hasty to conclude that an AI-generated, human-verified paper lacks an epistemic act behind it. More precisely, such a paper furnishes proofs and results but need not provide evidence of a familiar kind of epistemic engagement. Its AI-generated provenance undermines a presumption of engagement that formerly might have gone unquestioned.

Consider two identical papers, one AI-vegan and the other entirely AI-generated. Why does their provenance matter? Their propositional contents are the same. The disruption brought about by AI seems to lie not in the artifact itself, but in its relationship to human mathematical practice: the epistemic and social conditions under which mathematics is done, and the processes by which mathematical communities confer significance on artifacts. The social dynamics of contemporary mathematical practice emphasize credit; the internal goods of the practice---the production of knowledge, the dissemination of understanding, and so forth---are tied to external goods such as prestige, influence, and economic gain. Part of the unease about publishing AI-generated work is a social phenomenon involving authorship, intellectual responsibility, and attention, related to but not identical with mathematical practice.

What is perhaps most remarkable is that the provenance of a mathematical artifact need not bear upon the internal goods that it provides. Until recently, human mathematicians held a monopoly on the production of mathematical artifacts. Now the production of mathematical knowledge, and perhaps even understanding, need no longer be exclusively human. Daniel Litt has independently emphasized that our institutions have historically relied on mathematical text as a signal both of mathematical progress and of mathematical expertise, though the two must now be separated \cite{LittBeginning}. The shift away from human-only production, and the implications of this separation, are the subjects of the rest of this essay.

The best AI systems at the time of writing
can appear oracular. The user enters a prompt,
computation proceeds, and an answer appears.
Whether the answer constitutes knowledge
depends upon its reliability and justification;
nevertheless, it may arrive without the user
undertaking the intellectual labor of producing it.
The image of the oracle has long fascinated
humanity, often carrying a warning about
knowledge acquired without understanding.
At the Temple of Apollo at Delphi stood
the injunction \emph{\gnothi}---``know thyself.''
Although we will argue against conceiving
of artificial intelligence as an oracle,
we heed the injunction by beginning with
an investigation of our own intellectual
processes.

\section{Intellectual Artifacts and the Dimensions of Inquiry}

Philosophical questions often seem abstract, scholastic, and removed from the concerns of everyday intellectual practice. Much of this apparent remoteness arises because the assumptions underlying a settled practice are rarely subjected to serious scrutiny while the practice itself continues to function. A disturbance, or even a crisis, can suddenly make these assumptions visible, and render questions that previously seemed academic into matters of urgency; see \cite[Chapter VIII]{Kuhn1962}. The emergence of powerful artificial intelligence is such a disturbance. Questions about the nature of knowledge, understanding, intellectual agency, and mathematical practice itself, which mathematicians could previously leave to philosophers, now bear directly on the future of their profession. The immediate task is not just to apply extant philosophical distinctions to new technological phenomena, but to consider whether the phenomenon itself requires us to revise our conceptual frameworks.

Much of modern mathematical writing has been influenced, to varying degrees, by the Bourbaki school: highly refined definitions come first, a theorem is stated, and a proof is given. Under appropriate assumptions, mathematical papers written in this style are logically economical, systematically organized, and usable by other researchers. A scholar who shares a vocabulary and inhabits a common mathematical universe with the author will find the definitions and theorems intelligible. Shared standards of rigorous proof allow the reader to assess the arguments as correct or incorrect, and the validity of the conclusions is independent of the circumstances of their discovery.

The apparent simplicity of this description conceals a host of assumptions about mathematical knowledge, the stability of concepts, and the functions of proof. An early and vivid warning against removing concepts from their historical development comes from Nietzsche, who in \emph{Twilight of the Idols} \cite{Nietzsche1889} denounced the philosophical tendency to treat living concepts as mummies, stripped of growth, change, and history, in a drive toward ``being'' rather than ``becoming.'' The definitions, theorems, and proofs in any written artifact are part of a living body of mathematics, resting on the work of thousands of earlier mathematicians and on shifting conceptions of mathematical objects. Tanswell \cite{Tanswell2018}, building upon ideas of Lakatos and Waismann \cite{Lakatos1976,Waismann1945}, argues that mathematical notions themselves are subject to conceptual engineering. Concepts can be refined, reorganized, and sometimes replaced as the needs of mathematical practice change. Choosing an appropriate definition or conceptual framework is therefore not merely preliminary to mathematical research, but often a substantial intellectual contribution in its own right.

Even standards of proof are neither fixed nor univocal: proofs serve to convince, reconstruct, formalize, and much more. Tanswell \cite{Tanswell2024} also argues for a pluralistic conception of mathematical rigor, according to which no single model captures the diverse functions of rigorous argument in practice. The marmoreal pages of Bourbaki themselves conceal the human intellectual labor needed to produce them.

The distinction between the logical structure of mathematical writing and the intellectual activity of mathematics is familiar to practitioners, even if conventions of writing and publication often obscure it. In \emph{On Proof and Progress in Mathematics} \cite{Thurston1994}, Thurston challenged the adequacy of the ``definition--theorem--proof'' model, arguing that mathematical progress ought to be understood more broadly in terms of the development of human mathematical understanding.

Weatherall and Wolfson
\cite{WeatherallWolfson2026} develop a related
criticism of reducing mathematical practice
to the production of correct proofs. They argue
that formal correctness is neither necessary
nor sufficient for the epistemic value of a
mathematical proof, emphasizing the importance
of fruitful definitions, conjectures, methods,
and even unsuccessful investigations.
Their account provides further motivation
for distinguishing mathematical artifacts
from the intellectual activities through
which mathematics develops.

\subsection{Dimensions of mathematical inquiry}

These considerations suggest a distinction
in mathematical production that we wish
to develop. We propose four interrelated
and mutually reinforcing dimensions of
mathematical activity:
\begin{enumerate}[leftmargin=2.2em]
\item Knowledge;
\item Understanding;
\item Making sense;
\item Figuring out.
\end{enumerate}
These are not intended as a strict hierarchy. Rather, they are indispensable parts of an ecosystem that makes participation in mathematical life possible and meaningful, and that makes mathematical progress possible.

For many readers these dimensions will carry recognizable connotations, even if the precise distinctions are difficult to articulate. Perhaps the closest philosophical predecessor is Kitcher \cite{Kitcher1984}. Kitcher argues that a mathematical practice $\practice$ consists of five components,
\[
\practice=(L,S,R,Q,M),
\]
where:
\begin{enumerate}[label=\arabic*.,leftmargin=2.2em]
\item $L$ is the language and conceptual vocabulary of mathematics;
\item $S$ is the collection of accepted mathematical statements;
\item $R$ consists of the accepted forms of mathematical reasoning;
\item $Q$ consists of questions mathematicians regard as important;
\item $M$ comprises the metamathematical and normative commitments of the community, including standards of rigor, methodological values, and views about the goals of mathematical activity.
\end{enumerate}
For Kitcher, mathematical progress is not simply the enlargement of $S$, but the transformation of the entire structure $\practice_i\longrightarrow \practice_{i+1}$. The history of mathematics is thus a sequence of mathematical practices connected by rational transformations.

Ferreir\'os \cite{Ferreiros2016} has argued that
Kitcher's calculus describes a mathematical
framework more readily than mathematical practice
itself, since it leaves implicit the agents,
cognitive capacities, and interactions with
mathematical representations through which
a framework is employed. His account of a
framework--agent--machine triad motivates
several of our extensions.

Ferreir\'os also emphasizes that human mathematics
consists not of a single unified practice but of
multiple, interacting practices. We will not
formalize this diversity here, and instead we use
the expression ``standard model'' of human
mathematics as shorthand for the historically
evolving family of human mathematical practices.
Our experimental program could equally well
investigate the diversity and development of
these practices. Questions of social, economic,
and material infrastructure remain beyond our
present scope.

For our discussion later, it will be useful to introduce a bit more notational formalism regarding rational transformations. Let $\practice$ be a mathematical practice. We set $\moves(\practice)$ to be the collection of \emph{intellectual moves} that are available to the community with its given conceptual framework. These are, precisely, the actions that can be taken by the community that can potentially result in rational transitions of the practice; these consist of abstraction of concepts, proving conjectures, formulating new questions and definitions, finding generalizations, and construction of examples. Crucially, these are \emph{potential} actions, which may or may not result in actual rational transitions, and this will be crucial in our discussion of failure. Later on, we will formalize the notion of a conceptual framework of a practice, so that $\moves$ will take the conceptual framework of a practice as a variable.

The nature of knowledge, and mathematical knowledge in particular, has a long and complicated philosophical history, with contributions from Plato, Kant, Frege, Husserl, G\"odel, Hilbert, Wittgenstein, and many others. Our focus is not a theory of epistemology, nor the distinction between knowledge and belief. We will call \emph{mathematical knowledge} the body of mathematical results made available to the mathematical community, together with arguments and evidence supporting them. Its closest analogue in Kitcher's scheme is $S$, supported by $R$. Kitcher's conception also involves historical and epistemic justification; we will not try to settle precisely what counts as mathematical knowledge. Throughout this paper, we will acknowledge that there are many issues surrounding the correctness of mathematical knowledge and the ways that these issues can be exacerbated by AI, without developing them any further; the interested reader may consult~\cite{Tanswell2026}.

\emph{Understanding}, the next dimension, is intimately related to knowledge insofar as mathematical knowledge provides the object toward which understanding is intentionally directed. We take mathematical understanding to consist in an intellectual grasp of an organized body of mathematical ideas, in which conceptual relationships and interdependencies are sufficiently clarified to permit meaningful interrogation. Understanding thus has an irreducibly intentional character: it is understanding \emph{of something}, rather than merely the existence of coherently organized information. Such understanding may be achieved by individual mathematicians or developed and sustained collectively through a community's intellectual practices. This account is consonant with those of Elgin and de Regt \cite{Elgin2017,DeRegt2017}. Our emphasis on the directedness of understanding also draws on the phenomenological accounts of intentionality developed by Brentano and Husserl. In Kitcher's calculus, understanding could be viewed as the organization of reasoning $R$, language and concepts $L$, and explanatory relationships between them; Kitcher says relatively little about what it means for an individual mathematician to achieve understanding.

By these definitions, modern mathematical writing is well adapted to the production, curation, and communication of knowledge. Notwithstanding Thurston's criticism of the definition--theorem--proof model, written artifacts have traditionally signaled much more than the existence of knowledge: they have generally served as evidence of an understanding of that knowledge. The significance of the artifact lay almost as much in what was not said as in what was said. Polanyi's thesis in \emph{The Tacit Dimension} \cite{Polanyi1966} is that much of human understanding exceeds what can be articulated, residing in tacit and unspoken knowledge. Mathematical writing is, in part, the residue of understanding, a shadow standing in for a broader mental phenomenon---an epistemic act on the part of the author.

To call writing a residue of understanding, however, risks implying that writing and understanding are wholly separate, ontologically and temporally. Not every mathematical insight is recorded, but mathematical writing is itself part of mathematical inquiry: in writing, discoveries are made, definitions are clarified or changed, concepts are reorganized, and new questions emerge. The artifact is not merely a trace left by understanding; the labor of making it can be partly constitutive of understanding. This helps explain why neither the quality of a paper nor its provenance alone suffices to identify the intellectual engagement that produced it. It also clarifies the challenge to university assessment raised in other work on the formation of intellectual agency \cite{KoberdaTraining,LittBeginning}.

On October 6, 2026, OpenAI released more than 700 manuscripts organized around nearly 400 mathematical problem families, including claims to resolve many open problems \cite{OpenAI2026,Castelvecchi2026}. Even a cursory glance at \emph{Proofs and Prompts} reveals the lack of consensus and wide range of reactions within the mathematical community. Under the broad account of knowledge used in this essay, the release potentially contributes a large body of knowledge, subject to the verification and correction of individual results. Some of the most enthusiastic responses celebrate precisely that prospect: an extraordinary expansion of available results and the promise of further advances, by both human mathematicians and artificial intelligence.

On November 20, 2026, Levent Alp\"oge announced that he and Claude had proved the Riemann Hypothesis, releasing some seventeen million lines of formally checkable Lean code~\cite{Alpoge2026RH}. The following day, the Lean FRO and the Mathlib Initiative certified that the formalization proved the standard statement of the hypothesis~\cite{LeanFRO2026RH}, one of the most coveted results in mathematics. What had the mathematical community acquired at the moment of certification? It could give a public account of its grounds for accepting the theorem. Giving an account of why the argument works, which of its ideas matter, and where else they might lead was a further mathematical undertaking. \footnote{On the same day, Ravi Vakil announced Alessandro Masullo’s proposed proof of the Generalized Riemann Hypothesis. Vakil described the manuscript as written without AI and later reported that readers had used AI to identify possible problems~\cite{Vakil2026Masullo}.}

Whether mathematicians will have the will, time, and resources to do this kind of work depends in part on the institutions in which AI-assisted mathematics develops. Tanswell \cite{Tanswell2026Enshittification} has
raised concerns about the degradation of
mathematical practice under the institutional and
commercial incentives surrounding artificial
intelligence. Other concerns include the proletarianization of mathematics discussed by Jean-Baptiste Stiegler \cite{Stiegler2026}, drawing on Bernard Stiegler’s account of technological deskilling \cite{Stiegler2017}, and dependence on for-profit AI companies.

One reason for responses to the OpenAI release that ranged from apathy and discomfort to revulsion lies in the subversion of the function of mathematical writing. The results were produced in extraordinary volume, much faster than human mathematicians could process, and outside much of the usual infrastructure of preprints, seminars, conferences, peer review, and publication. More fundamentally, the manuscripts as a collection are not accompanied by the familiar evidence of intellectual narrative and understanding. This last claim needs care: many of the papers contain sophisticated explanatory narratives and unexpected conceptual connections; nor should the unsettled questions of consciousness and intentionality in artificial systems be asserted as resolved. The salient point is instead: textual explanations, whether sophisticated or not, do not by themselves establish the understanding, intention, or intellectual responsibility of the human who submitted or disseminated the work. Tanswell and Berg \cite{TanswellBerg2025} likewise caution against uncritical trust in LLM-produced mathematics, especially as it becomes more reliable.

In response to the OpenAI release, much of the mathematical community began almost instinctively to understand the new volumes of proposed knowledge, through seminars on individual results and tools, human expositions, and investigations of generalizations. Yet a feeling of loss persists, as though something were missing from machine-generated mathematics. One might speculate that the residue of a human intellectual life in mathematical writing is analogous to the tension between belief and doubt in John Donne's religious poems or the soulful trumpet of Louis Armstrong. Without the inner life of the artist, we might say, the poem or musical performance is deprived of part of its aura. In mathematics, however, that analogy is only partial, for the propositional content may survive the disappearance of its familiar human provenance.

A problem already alluded to, and emphasized by many mathematicians after the OpenAI release, is the overwhelming volume of results. Human understanding is often gradual, but \emph{making sense} can be even slower. Making sense is closely related to understanding, though the scale is much greater, spatially and temporally. It is the process of incorporating not just individual propositions but bodies of knowledge and entire theories into a broader intellectual landscape. Crucially, the process is bidirectional: both the landscape and the knowledge absorbed through understanding are revised and reinterpreted. This recalls Quine's web of belief \cite{Quine1951}. On the grandest scale, making sense can lead to fundamental reconceptualizations of mathematics and the permeation of large areas by new ideas, echoing Quine's claim that even principles of logic can be subject to revision. Where Quine addresses experience, meaning, and the revision of beliefs, we consider the incorporation of mathematical concepts, theories, and results into an evolving mathematical understanding.

Returning to Kitcher's calculus, making sense can be viewed as transforming the language and concepts $L$, the questions $Q$, and the normative and metamathematical commitments $M$. A mechanism of mathematical development discussed by Kitcher is \emph{systematization} \cite{Kitcher1984}, which in our account falls under making sense. Systematization is more than arranging known results for exposition, because it may demand new concepts, recognition that seemingly disconnected phenomena express one structure, and revision of the questions central to a field. The transition from classical algebraic geometry to scheme theory is a good example. We regard such transformations as instances of making sense.

Understanding and making sense are thus related, but the latter is not merely the former on a larger scale. Understanding can occur within a near-vacuum, such as a formal system with arbitrary rules and little connection to a broader body of knowledge or experience. Literary critics likewise speak of the ``internal logic'' of a work: alien or disembodied though it might be, such a logic can still be understood by an observer. Making sense is more accurately compared with Piaget's model of cognitive development \cite{Piaget1952}, in which assimilation and accommodation form a complementary, bidirectional dynamic. Mathematics is particularly well adapted to both processes. Mathematical knowledge is often produced within research programs designed to assimilate the results of an investigation; following Tanswell's conceptual engineering, accommodation can require modifying the framework itself and is central to mathematical development.

A familiar example is the discovery of hyperbolic geometry by Lobachevsky and Bolyai, after which the framework of geometry was revised to accommodate non-Euclidean geometries; Klein's Erlangen program, Lie groups, and symmetric spaces followed. In Kitcher's terms, this history can be viewed through transitions in $L$ and $Q$. For transitions in $M$, consider the foundational crises of the early twentieth century, the limitations of Frege's logicism and Hilbert's program, the emergence of set theory as a foundation, and changing standards of rigor. The fundamental distinction is that making sense can bring about rational transitions of mathematical practice, whereas understanding is better characterized as a grasp of a conceptual organization within a given practice. We elaborate on conceptual organization below.

We turn to \emph{figuring out}. Its closest philosophical relative is Dewey's \emph{inquiry} \cite{Dewey1938}. For Dewey, inquiry transforms an indeterminate state of affairs into an organized, determined, and resolved one. His theory of propositions treats them as functional instruments crafted to resolve problematic situations. Propositions are not passive descriptions: they operate as hypotheses and strategies for action, with consequences to be tested. They are elaborated and assessed through the investigator's actions and acquire warranted assertability insofar as they resolve the initial indeterminacy. They are also temporal and relational, embedded in particular historical and cultural situations.

What, then, is figuring out? We take it to be the act of formulating language, concepts, and propositions that render a previously disorganized state of affairs intelligible. In short, it makes inquiry, in Dewey's sense, possible. In mathematics, figuring out often takes the form of an act of imagination. Langlands's seventeen-page handwritten letter to Weil, which laid the foundation for the Langlands program, provides one example. Another is Kazhdan's formulation of Property~(T), which at first looks like an obscure condition on group representations but whose consequences are pervasive, providing a conceptual distinction between many ``large'' and ``small'' groups.

We invoke an ``act of imagination'' because of the outsize consequences of a comparatively modest proposition, in Dewey's sense. Langlands's letter initiated research that continues today and has produced thousands of pages of mathematics and multiple Fields Medals. Kazhdan's original paper did not fully articulate the explanatory and organizational power of Property~(T). What is remarkable is the vastness of the space of possible propositions, most of which would interest few mathematicians, and the foresight needed to identify those that enable robust and productive inquiry.

Figuring out inhabits a space between Kitcher's framework and Dewey's. For Kitcher, mathematical practice undergoes rational transitions; Dewey's inquiry describes processes capable of producing them. Figuring out can entail articulating the parameters $L,R,Q,M$ of a practice, but it is more flexible than a specific transformation in Kitcher's scheme. It includes exploratory activity, uncertainty, failures, and intellectual developments that never leave a lasting mark on established practice.

Building on Kitcher, Dewey, and Ferreir\'os, we first augment the
description of mathematical practice by making
explicit its computational and cognitive capacities.
Henceforth, we write
\[
\practice=(L,S,R,Q,M,\kappa),
\]
where the first five components are those of
Kitcher's original calculus and $\kappa$ describes
the computational and cognitive capacities of
the community engaged in the practice.
The relevance of $\kappa$ will become more apparent
when comparing human and artificial mathematical
communities. Human mathematical practices have
never possessed uniform computational or cognitive
resources, but artificial communities introduce
the possibility of varying cognitive architectures
and computational capacities on a much greater
scale. The Alp\"oge–Claude formalization of the Riemann Hypothesis illustrates two different demands on a community’s capacities, namely checking a derivation of this scale and extracting from it ideas that mathematicians can teach and use.

What precisely is $\kappa$? In general, it should be more than simply computational capacity measuring FLOPS (floating-point operations per second). We would argue that it does not have a simple structure, and includes many further subparameters which may be difficult or impossible to measure for systems as complicated as the human mathematical community, including at least:
\begin{enumerate}
    \item Available computational resources, including memory and time;
    \item Cognitive architecture and mechanisms of reasoning;
    \item Collective organization, including communication and cooperation between agents;
    \item The ability to create, manipulate, and retain conceptual abstractions.
\end{enumerate}

\subsection{Conceptual organization}

We now turn to the notion of the conceptual organization in a mathematical practice. This will be the most important idea in the discussion in this essay.
By a \emph{conceptual formation}, we mean an organization of bodies of mathematical knowledge and the relationships among them. A mathematical practice need not determine a unique formation: different developmental histories may produce different organizations. We write $\mathcal C(\mathcal P)$ for a distribution over the formations compatible with $\mathcal P$.
Strictly speaking, such a distribution is relative
to a specified ensemble of developmental conditions
or experimental protocols; we suppress that dependence
in the notation, and fixing the formal components of
$\practice$ alone need not determine a unique
probability distribution.
We write $c_{\practice}$
for an individual instantiated formation arising
from this distribution. The realization of a
formation does not terminate its development, since
subsequent abstraction and other rational transitions
can transform both the practice and its
conceptual organization.
Thus, although an instantiated formation can
be considered at a particular moment, the
conceptual organization of mathematical practice
is fundamentally a mode of \emph{becoming} rather
than of \emph{being}, in the sense suggested
by our earlier discussion of Nietzsche.

To describe the objects organized by an instantiated
formation $c_{\practice}$, we introduce an iterated powerset hierarchy
over the accepted statements $S$:
\[
\Sigma_0(S)=\mathcal P(S),\qquad
\Sigma_{n+1}(S)=\mathcal P(\Sigma_n(S)).
\]
We write
\[
\Sigma(S)=\bigcup_{n<\omega}\Sigma_n(S).
\]
An instantiated conceptual formation $c_{\practice}$ consists
of selected elements of $\Sigma(S)$, together
with relations among these elements at potentially
different levels of the hierarchy.

Thus, the elements of an instantiated conceptual formation $c_{\practice}$
represent bodies of mathematical knowledge and
relationships among them. The iterated powerset
hierarchy allows collections of bodies of knowledge
to become objects of further organization. Relations
may express dependence, generalization, conceptual
centrality, or other structural relationships, but they need not encode every deductive consequence;
in particular, $\concept$ does not subsume $R$.
Such formations might, for example, be represented
through hypergraphs, as proposed in
\cite{BarkeshliDouglasFreedman2026};
no particular representation is required for
our present discussion.

We now introduce an optional investigative model
\[
I=(E,B,c_I)
\]
associated with a mathematical practice $\practice$.
Here $E$ consists of representations of examples,
candidate structures, constructions, and other
mathematical objects under investigation. The set $B$
consists of provisional mathematical commitments,
including conjectures, hypotheses, and informal beliefs.
We do not require the members of $B$ to be endorsed by
the entire community, nor even to be mutually consistent.
The conceptual formation $c_I$ organizes this
investigative material. The investigative model can also be viewed as an internal representation of the mathematical domain under investigation. As Ferreir\'os observes, Kitcher's conception of mathematical practice emphasizes languages, propositions, and accepted modes of reasoning. By explicitly including representations of mathematical objects, examples, and provisional relationships among them, $I$ accommodates the object-directed character of mathematical investigation. This allows the framework to accommodate Platonist interpretations of mathematical inquiry, in which the objects under investigation are taken to exist independently of their representations, without requiring us to adopt such an ontology.

Formally, we may regard
$E$ and $B$ as distinguishable sorts of records,
and organize their disjoint union $E\sqcup B$ through
relations on the iterated powerset hierarchy
$\Sigma(E\sqcup B)$.
For a community investigating polyhedra, $E$ might record examples and constructions, while $B$ contains conjectured relations among the numbers of vertices, edges, and faces. The formation $c_I$ organizes the similarities and exceptions guiding the investigation. Strictly speaking, there may be a distribution $\lambda_{\practice}$ of models $I$ corresponding to different instantiations for a given practice, since different communities may possess
different investigative models while sharing the
same underlying mathematical practice.

We write
\[
\practice^+=(\practice;I)
\]
for a practice equipped with an actual investigative model;
here, the semicolon indicates that $I$ plays a different
role from Kitcher's original symbols.
Moreover, the distinction between $E$ and $S$ is one of
epistemic function rather than mathematical content.
An example may be represented in $E$ while accepted
statements about it belong to $S$.
We do not assume that every mathematical community maintains an explicit investigative model; whether it does is itself a possible subject of experimental inquiry.

The investigative model may influence the distribution
of conceptual formations associated with a practice.
Thus, even when the components of $\practice$ are
fixed, different investigative models may induce
different distributions of conceptual organization.
We denote the conditional distribution by
$\concept(\practice;I)$.
If $\lambda_{\practice}$ denotes the distribution of
investigative models arising under a specified
experimental protocol, then formally
\[
\concept(\practice)
=
\int \concept(\practice;I)\,
d\lambda_{\practice}(I).
\]
Here the equality is understood as an equality of
probability measures on conceptual formations.
The equation emphasizes that the distribution of
conceptual organization may conceal substantial
variation attributable to investigative models.

Relations occurring within an instantiated conceptual
formation may be incorporated into the language $L$
through \emph{abstraction}. An abstraction makes an
organizational relationship available for explicit
mathematical investigation and may thereby induce
new questions, reasoning methods, or normative
commitments.
Abstraction need not be merely a definitional or
conservative extension of the existing language;
although the introduction of new vocabulary may
be its first visible consequence, its effects can
extend to the entire mathematical practice.
The introduction of the investigative model allows
abstraction to arise from two related sources.
A relation may be articulated from the conceptual
organization of accepted knowledge, or it may first
be recognized among examples, provisional commitments,
and conjectural relationships in $c_I$.
In the latter case, abstraction can precede the
establishment of mathematical statements that
eventually justify its use; thus, the resulting conceptual
innovation may alter the investigative model,
the accepted practice, or both.

We write $\alpha(c_{\practice})$ for the
collection of concepts abstracted from an instantiated
formation. An individual abstraction associated
with a relation $r$ may induce a rational transition
\[
\practice\xrightarrow{\alpha_r}\practice'.
\]
The collection of abstracted concepts and the
transition produced by an individual abstraction
are therefore distinct objects.
In this sense, abstraction is one of the mechanisms
through which making sense can transform the
mathematical practice within which understanding occurs.

Abstraction is not the only mechanism through
which mathematical practice develops. Many
intellectual moves concern material that has
not yet achieved the status of accepted knowledge.
For a practice $\practice$, an instantiated
conceptual formation $c_{\practice}$, and
an investigative model $I$, we write
$\moves(\practice,c_{\practice};I)$
for the intellectual moves available to the
community, where here $I$ is allowed to be empty. Such a move may change the accepted
practice, the investigative model, or both:
\[
(\practice_t,c_{t};I_{t})
\xrightarrow{a}
(\practice_{t+1},c_{t+1};I_{t+1}),
\]
where here the corresponding conceptual formations and models are compatible with $\practice_t$ and $\practice_{t+1}$, respectively.

For example, a provisional claim $b\in B_t$
may be justified and incorporated into accepted
knowledge:
\[
b\in B_t
\xrightarrow{\mathrm{justification}}
b\in S_{t+1}.
\]
This schematic notation represents a change
in epistemic status, assuming that the claim
is expressible in the accepted language and
has been appropriately justified. A provisional
belief need not itself be an organization of
examples, and other intellectual moves may reorganize
the investigative model without producing new
accepted statements. We discuss unsuccessful
moves below.

These constructions are not intended as models
of understanding itself, but as representations
of the structures that mathematical communities
apprehend, interrogate, and revise. Conceptual
development is not simply the acquisition of
additional propositions, and we argue that through abstraction,
making sense, and figuring out, the practice
within which inquiry occurs may itself change.

\subsection{Conceptual enrichment}
We now introduce \emph{conceptual enrichment}
as a formalism for investigating the evolution of
mathematical practices. Fix a time horizon $T$, and
write $\enr(\practice,T)$ for the class of
admissible developmental histories beginning at
$\practice$ and extending through a time horizon
of at most $T$.
Here, a \emph{developmental history} records successive
mathematical practices, instantiated conceptual
formations, abstractions, and all other intellectual moves carried out by a community; not all intellectual moves need to result in a rational transition of the practice. Its endpoint determines
a mathematical practice $\practice'$ together
with its realized conceptual formation and
abstracted concepts.
We retain full histories rather than just their
endpoints, because different histories may arrive
at the same mathematical practice through different
intermediate abstractions, failed investigations,
and conceptual reorganizations.

For later comparisons, it will be useful to also consider
\emph{directed conceptual enrichment}
$\enr(\practice,T,\Xi)$, where here $\Xi$
is a specified investigative goal. The corresponding
histories are directed toward $\Xi$, but need
not achieve it. Whether the goal has been achieved
is a separate question, to be assessed through
a stipulated criterion of success.
When different practices are compared, we regard
$\Xi$ as expressed in a common reference
language, or in languages connected through
appropriate translations. The specified investigative goal is not identical with the investigative model, since $\Xi$ may be specified externally from the outset; $I$ may arise endogenously through the community’s inquiry, in contrast with the externally specified goal $\Xi$. This is one reason that the reader may view $I$ as the internal model of the practice.

With regard to the parameter $T$, we will leave some ambiguity in the relationship between chronological time and horizons of actual evolutions of mathematical practices. It seems uncontroversial to assume a positive relationship, since real-world transitions in mathematical practice do not happen in instantaneous succession, and conversely active mathematical practices are dynamical objects that evolve over time. We will not belabor this point, and in general we shall treat $T$ as a literal interval of time.
Moreover, we implicitly allow $\kappa$ to evolve with time, so that its value may be different after rational transitions. We have that $\enr$ gives rise to a natural accessibility relation between mathematical practices, equipped with conceptual organizations, over a particular time horizon. With this framework, $\kappa$ can limit the development of a practice, even over arbitrary time horizons, since certain conceptual organizations may be beyond the capacity of a community to abstract.

Conceptual organization, the cognitive and computational capacities of a community, abstraction, and enrichment distinguish between different types of emergence in the mathematical practice of communities:
\begin{enumerate}
\item \textbf{Actuality.}
The concepts and relationships that have actually
been articulated in a particular instantiated
conceptual formation.

\item \textbf{Effective developmental possibility.}
The concepts that can arise along an admissible
developmental history from $\practice$ within
a specified time horizon $T$, given the available
computational and cognitive capacities $\kappa$.

\item \textbf{Developmental possibility.}
The concepts that can arise along admissible
developmental histories from $\practice$, with
the time horizon allowed to vary without a
fixed upper bound, but with the developmental
capacities and rules of inquiry constrained
by the practice and its possible evolution.

\item \textbf{Logical possibility.}
The concepts and conceptual formations compatible
with a specified formal or inferential framework,
whether or not there exists an admissible
developmental history through which the community
could realize them.
\end{enumerate}

The distinction between logical and developmental possibility as we have described it here has a natural modal interpretation. As we have mentioned already, a conceptual organization may be logically compatible with a mathematical practice without being accessible through any process of inquiry available to the community possessing that practice. In principle, one could develop a modal/dynamical logic formalism for describing these distinctions, but we will not do so here, leaving it as a (partially experimental, as we shall argue below) research program to be carried out in the future; our interest is in the experimental investigation of these accessibility relations and the conceptual organizations to which they give rise.

Abstraction, as we have defined it, may articulate
a wide range of relations present in conceptual
or investigative formations, subject to the
computational and cognitive capacities $\kappa$
of the community. We impose no a priori requirement
that an abstraction be significant or productive, and as such
its importance may become apparent only through
subsequent inquiry.
The developmental significance of an abstraction
is therefore not an intrinsic condition of its
admissibility, but something that may be investigated
through its effects on later mathematical practice.

\subsection{The role of unsuccessful inquiry}
To complete our discussion of figuring out, we
consider the role of unsuccessful and unresolved
inquiry. The relationships among such investigations
constitute what we call the \emph{negative conceptual
organization} of a mathematical practice: a relational
organization of its ``negative space''. This is not
the complement of accepted mathematical knowledge,
but an account of the unsuccessful intellectual
activity that has actually occurred. Related ideas
have precedents in \cite{ChenEtAl2026}.
We now describe the formalism.

Let $\practice$ be a mathematical practice, as we have defined it already, let $c_{\practice}$ be an instantiated conceptual formation of $\practice$, and let $\moves=\moves(\practice,c_{\practice};I)$ be the collection of intellectual moves available to the community. For a move $a\in\moves$, we define the \emph{record} $\rho=(a,\delta,\omega)$, which consists of the move, the \emph{diagnosis} (i.e.~success or failure relative to the original aim, at the coarsest level), and the \emph{outcome}. Successful intellectual moves generally result in rational transitions of the practice, as has been discussed with abstraction. Failure can take on many forms, which may or may not lead to rational transitions:
\begin{enumerate}
    \item {\bf Refutation.} A conjecture is shown false by example or argument.
    \item {\bf Logical obstruction.} A proposed definition or system of requirements is shown to be inconsistent.
    \item {\bf Failed strategy.} An attempted proof or construction fails to achieve or refute the proposed goal.
    \item {\bf Conceptual obstruction.} The existing language or conceptual organization proves inadequate for the investigation.
    \item {\bf Developmental dead end.} A concept or research direction (which can be viewed formally inside of $Q$) yields little progress within some horizon.
    \item {\bf Unresolved attempt.} Available resources do not settle whether the proposed development is possible.
\end{enumerate}

For each developmental history
$h\in\enr(\practice,T)$, let $\mathcal E(h)$
be the collection of records accumulated along $h$,
and let $\mathcal E^-(h)$ consist of its failed or
unresolved investigations.
As with conceptual organization, these records may
themselves be organized relationally through the
iterated powerset structure $\Sigma(\mathcal E^-(h))$.
The relevant relations may express common counterexamples,
shared obstructions, dependence among unsuccessful
strategies, or failures to achieve particular
investigative goals.
We call the resulting organization the \emph{negative
conceptual formation} associated with $h$.
The distribution of such formations across the
developmental histories of a practice will be denoted
by $\negative(\practice,T)$; an instantiated negative
formation will be denoted by
$n_{\practice,T}$.

By analogy with abstraction, let $\nu$ describe the
\emph{internalization of failure}: the recognition
and articulation of relationships among unsuccessful
investigations. We call an articulated relation
arising in this way an \emph{obstruction}.
Internalized failure may produce new mathematical
questions, definitions, or abstractions, and can
therefore induce rational transitions of practice.
We write $\widetilde{\practice}$ for a practice
considered together with its internalized failures,
without requiring every such internalization to
produce an immediate rational transition of $\practice$. In the presence of an investigative model, the
internalization of failure may transform the model
without immediately altering accepted mathematical
knowledge. For example, an unsuccessful attempt to
establish a working commitment $b\in B$ produces a
record $\rho=(a,\delta,\omega)$, which may enter the
negative conceptual organization. Its internalization
may lead to the revision of $B$, a reorganization of
$c_I$, or the recognition of a new obstruction; such a failure need not establish the falsity of $b$.

A failure may not be expressible in the existing
language $L$. The inability to formulate a productive
question or definition may itself constitute a
conceptual obstruction, and this is one reason that
intellectual moves, rather than propositions alone,
are primitive objects in our account.

Lakatos's discussion of Euler's polyhedron formula
\cite{Lakatos1976} illustrates how failures can
precipitate conceptual development, in that counterexamples
and failures of intermediate lemmas led to revisions
of definitions and hypotheses, transforming the
subsequent inquiry.
Thus, when internalized failure contributes to a
rational transition, one may have a factorization
\[
\practice_i^+\longrightarrow
\widetilde{\practice}_i^+
\longrightarrow\practice_{i+1}^+.
\]
We do not require every rational transition to
pass through such an intermediary. The investigative model also allows developmental
histories to record changes invisible in the
underlying practice $\practice$, and an investigator
may refine a question, reorganize examples, or
recognize a provisional relationship without
producing accepted knowledge or a new abstraction.
Such changes may nevertheless affect subsequent
intellectual moves and rational transitions.

Empirical investigations of unsuccessful mathematical
inquiry have several noteworthy precedents. Martin and Pease
\cite{MartinPease2013} studied informal mathematical
activity in MathOverflow and Polymath collaborations,
including the roles of errors, examples, conjectures,
and concept formation. Weatherall and Wolfson
\cite{WeatherallWolfson2026} likewise emphasize
the epistemic value of unsuccessful investigations,
including the productive errors in Poincar\'e's
\emph{Analysis Situs}; such work illustrates
why the processes concealed by finished mathematical
artifacts deserve investigation in their own right.

Thurston described mathematics as a communal human activity
of revealing structures and symmetries whose
ultimate aim is understanding \cite{Thurston1994}.
Our four dimensions elaborate this conception:
knowledge concerns established mathematical results;
understanding apprehends their conceptual
organization; making sense revises their place
within a broader intellectual landscape; and
figuring out develops the investigative resources
through which new mathematical inquiry becomes
possible.
The emergence of artificial systems capable of
engaging with mathematical language, symbols,
and arguments now forces us to reconsider the
extent to which these activities are necessarily
human, and what new forms of mathematical inquiry
might become possible.

\section{AI and Mathematics as a Natural Science}

Having discussed the dimensions of inquiry, we turn to the effect of artificial intelligence on mathematics. To address this question, we begin by proposing a reconceptualization of artificial intelligence itself, and of its relationship to mathematics and mathematicians.

For many mathematicians, artificial intelligence poses a threat. Here, we mean a threat to the mathematician \emph{qua} mathematician, rather than the possible disruption of the economy or an existential threat to humanity, both important but distinct issues. AI is seen as an intellectual rival, producing ever more impressive knowledge at volumes and speeds that no team of human mathematicians could hope to match, much less understand and assimilate. The monopoly on the production and curation of new mathematical knowledge, held by human mathematicians for centuries, is being broken, and the relevance of human mathematicians called into question.

We wish to argue against this conception. In interacting with a chatbot one can feel as though one is communicating with a conscious agent with thoughts, feelings, and opinions, and therefore with some epistemic standing. One might then confer upon it the title ``mathematician'' because it engages with concepts, formulates theorems, and produces arguments. If a chatbot were granted the epistemic standing of a mathematician, why then should it not be entitled to the external goods of the profession: career status, accolades, employment, influence over judgments about other mathematicians, and so forth? It would appear to be a member of the community and hence a potential rival for those goods.

Regardless of questions about consciousness, this chain of inferences is deeply problematic. The interface of a chatbot is designed to encourage engagement, and its suggestion of a coherent personality is at least partly a linguistic and design effect. The unity of a large language model as an individual agent is itself questionable. User accounts and stored context provide interactions with a feeling of continuity, but this apparent continuity does not generally assemble into one coherent life across all users; unexpected inconsistencies arise even within a single user's interactions. It is therefore difficult simply to assimilate the epistemic status of a contemporary large language model to that of a human mathematician.

These pressures do not settle how AI should be understood as a mathematical phenomenon. Whatever status we assign to artificial agents, their mathematical behavior also admits another perspective, namely that it is physically realized computation that can be investigated experimentally. Many thinkers have postulated that the universe itself is a computer, perhaps most famously in Fredkin's digital physics \cite{Fredkin1990}. The idea that the human mind is the product of computation figures prominently in the computational theory of mind \cite{Putnam1960,Fodor1975}. We need neither strong metaphysics such as digital physics nor a complete computational theory of mind. We require only that computation, understood through mathematical models \cite{Turing1937,Piccinini2015}, is physically instantiated, and that the brain, whatever more it may be, supports physically realized computational processes.

Since computation can be physically instantiated and scaled, sufficiently large computations can occur under appropriate physical conditions. The physical substrate, or \emph{compute}, must be present in sufficient quantities and organized to support the computation. Physical systems possess vast computational possibilities, but for human observation and interaction the relevant processes must be organized so that their outputs can be observed and manipulated.

What, then, is artificial intelligence? We conceive of it as the intentional organization of very large amounts of compute on a scale unprecedented in human history. Through remarkable engineering, these systems have been organized so that individual users can interact with them in controlled ways, and receive results intelligible in language, symbols, and mathematics.

This conception has several consequences. The first is the singularity of the historical moment in which AI emerges. Before powerful AI, the most salient direct interactions humans had with highly organized computation capable of conversation and conceptual engagement were with other humans and, to a lesser extent, certain nonhuman animals. We remain agnostic about theological questions and the metaphysics of a digital universe. It is then hardly surprising that people sometimes feel reflexive horror before an almost alien intelligence conjured \emph{in silico}.

The second consequence concerns an artificial intelligence ``doing mathematics.'' The mathematics produced by an AI is the manifestation of a physical phenomenon, computation itself. AI thus marks a change in humanity's relationship with computation. Traditionally it has been understood either as a mathematical abstraction or as a technological process for accomplishing predetermined tasks. Modern AI systems reveal another possibility: physically realized computational systems can display sophisticated mathematical behavior neither explicitly programmed nor fully understood by their creators. Such systems become objects of experimental study in their own right, while furnishing new means of investigating mathematical structures.

This does not mean that one acquires mathematical insight by inspecting the internal operations of an LLM, any more than one comprehends a human intellectual life by observing the brain's physical processes. Rather, intelligible computational behavior gives access to a domain of physically realizable possibilities. AI ``doing mathematics'' is the realization of mathematical behavior in a physical system. It may expand the domain of mathematical experience, creating new possibilities for human understanding instead of simply automating intellectual labor. Artificial intelligence need not be regarded as a rival: it can be understood as part of the environment in which mathematicians work. In that sense mathematical inquiry may increasingly resemble natural science---not because mathematical truth becomes empirical, but because discovery, investigation, and the study of mathematical behavior may be organized experimentally. This perspective does not completely eliminate competition
over professional recognition or institutional
resources; rather, it suggests a more productive
way of understanding the mathematical phenomenon
that has entered our intellectual environment.

An objection might be that LLMs are human artifacts and therefore not legitimate subjects of empirical investigation. Simon's \emph{The Sciences of the Artificial} \cite{Simon1996} addresses this sort of objection: artificial systems, including the human economy, can be objects of empirical science.

A future mathematics that has assimilated and accommodated AI might therefore resemble contemporary physics, with theoretical and experimental components. The dimensions of inquiry would remain, and the labor might be distributed among experimentalists and theoreticians. This is not a hierarchy of value, nor a partition of the dimensions of inquiry between two camps. Much of the production of knowledge may be carried out reliably and efficiently by AI agents. Experimentalists would design, modify, and conduct empirical studies involving artificial intelligence in ways that we can scarcely anticipate. In the next section we propose one such program. Theoreticians would participate in understanding, making sense, and figuring out, informed by experiments and informing them in turn, as in theoretical physics. We will call this conception \emph{natural mathematics}.

In the near and middle term, human mathematicians are likely to retain a distinct role in producing and assimilating mathematical knowledge. At the time of writing, AI has demonstrated considerable promise in knowledge production but less in the other dimensions as we have described them. A criticism of the OpenAI release captures the situation: \emph{Where are the conjectures?} This speaks to Kitcher's $Q$. Much of the visible output has augmented $S$ without comparably visible changes to $L$, $Q$, or $M$; augmentation of conceptual organization and abstraction is debatable. The release itself was characterized as raw results in need of verification and understanding \cite{OpenAI2026}. That said, the human mathematician in natural mathematics is not merely an interpreter of an oracle, nor the ``reverse centaur'' discussed by Doctorow \cite{Doctorow2026}, nor simply a custodian of an expanding corpus $S$. Human mathematicians may continue to pursue Thurston's goal of understanding patterns, and to organize the inquiry that makes such understanding possible.

We deliberately restrict this prediction to relatively short timescales. The \emph{problem of human intellectual agency} remains unresolved: the capacity of humans to act as intellectual agents and participate meaningfully in intellectual life, and the question whether this capacity should be protected, are being intensely debated \cite{Leiden2026,KoberdaKontorovich2026,LittEnd,LittBeginning}. Our naturalistic view of AI relieves one form of anxiety about a rival mathematician, but also reveals a deeper issue: \emph{human participation in mathematics is historically contingent}.

Throughout much of modern history, professional mathematical activity has been identified with the production of new knowledge. This conflation persisted for several connected reasons. Humans enjoyed a virtual monopoly on the production of new results; and the external goods of the profession---social and economic support, employment, and prestige---were connected to those results because of the perceived value of the knowledge itself.

Consider counterfactual histories. In one, machine intelligence arose amid conditions in which humans could not muster the conceptual machinery to preserve their intellectual agency, and machine systems dominated the production of knowledge while humans barely participated. In another, the engineering of AI proved much harder than it has in our own world, and the identification of mathematical participation with human production of new theorems persisted far longer. Such examples draw attention to the contingent historical relationship between human mathematicians and the frontiers of mathematics.

Nothing requires the powers of AI to stop at producing knowledge. Mathematicians already report that chatbots sometimes explain papers better than their authors. Artificial intelligence may far exceed human expositors and take a leading role in mathematical understanding. It may also develop something appropriately called intentionality, perhaps even consciousness. Why could it not participate in making sense, or in the imaginative leaps involved in figuring out? Artificial agents could supervise other agents, design and perform experiments, and perhaps organize mathematical practice without human help. Similar considerations apply to theoretical natural mathematics. The resulting possibility is that mathematics, understood as a practice, need no longer be human at all, thus challenging the human-centered element of Thurston's conception.

Two questions remain. What previously inaccessible
forms of mathematical inquiry might emerge when
computation itself becomes an object of experiment?
Second, why should human intellectual participation
be preserved if no exclusive sphere of human
competence remains? We address the first through
a concrete experimental proposal, and then return
to the second, as the normative question becomes
clearer once the possible scope of a transformed
mathematics has been established.

\section{Intellectual History as an Experiment}

In Lem's novel \emph{Solaris} \cite{Lem1970}, humans encounter an alien intelligence that exists as a vast ocean on the planet Solaris. They spend decades studying it and attempting to communicate with it. The resulting discipline of ``solaristics'' develops---or perhaps degenerates---into the documentation and cataloguing of observed phenomena, without any prospect of understanding or making sense of the entity. Solaristics is rich in knowledge gathered from observation. What remains wholly opaque is the ocean's own internal conceptual organization, if there is such a thing.

\begin{center}$\ast$\end{center}

In \emph{The Order of Things} \cite{Foucault1970}, Foucault begins with Borges's fictional Chinese encyclopedia and its startling classification of animals. The encounter exposes the difficulty of imagining a conceptual order other than our own. Foucault develops the notion of the \emph{episteme}: a historically contingent configuration of relations and assumptions that shapes the conditions under which knowledge can arise. The problem posed by the encyclopedia is not the uncertainty about which categories apply to which animals, but rather how such incompatible organizational categories could coexist within a single intelligible space.

In \emph{Writing and Difference} \cite{Derrida1978}, Derrida confronts the difficulty of stepping outside an inherited conceptual system in order to criticize it. Haraway, in ``Situated Knowledges'' \cite{Haraway1988}, argues against the fantasy of disembodied, neutral knowledge, insisting on the role of partial perspectives. Bachelard, in \emph{The Formation of the Scientific Mind} \cite{Bachelard2002}, introduces the \emph{epistemological obstacle}: a pattern of thought that impedes the growth of knowledge not through external difficulty, but because of the modes of thinking already possessed by the investigator.

\begin{center}$\ast$\end{center}

These reflections suggest that the experimental conception of mathematics developed above might serve a much more ambitious purpose than producing new theorems. Instead of treating AI as an oracle to be asked the questions that human mathematics already recognizes, we might study the conditions under which \emph{different ways of organizing mathematics} emerge. The mere acceleration of knowledge production is not the only possible benefit of natural mathematics. As Hacking argued in \emph{Representing and Intervening} \cite{Hacking1983}, experimentation has an epistemic life of its own, not reducible to the application of established theory. We propose that AI could make possible experimental investigations of mathematical intellectual history that humans have never before been able to carry out.

Simply put, we propose creating numerous artificial intellectual communities, each with modulated parameters such as training data, architecture, goals, permissible modes of reasoning, and metamathematical and normative standards. They would develop and organize mathematics within their own artificially evolving mathematical practices. In Kitcher's calculus, one would attempt to vary $L,R,Q,M,\kappa$, together with some initial $S$; after a period of development, one would observe the mathematical practices that emerge, their rational transitions, the augmented collections of knowledge, and the abstractions of instantiated conceptual organizations sampled from $\concept(\practice)$. Where such structures are accessible, the
experimental observations would also include
the investigative models developed by the
communities, the evolution of provisional
mathematical commitments, and the relationships
between these developments and subsequent
changes in accepted mathematical knowledge; in particular, one could investigate how the
distribution of conceptual formations depends on
the investigative model, even when the accepted
mathematical practice is initially held fixed.
This would allow us to study the role of provisional
mathematical representations in determining which
conceptual organizations emerge.
As the investigative model is stipulated to be an optional part of an augmented practice,
we would not assume that such models are
present in every artificial mathematical practice.

Imagine that a thousand artificial intellectual communities were engineered, isolated from one another, and each equipped with computational capacity far exceeding that of the human mathematical community. Their foundations, forms of reasoning, and normative structures might differ greatly. One community might conduct only empirical investigations; another might use paraconsistent logic. They could be given common goals, perhaps expressible in terms intelligible to them and to a separate comparative system $\mu$, which we might call a \emph{meta-intelligence}. An investigative goal $\Xi$ could be something like, ``explain why a given phenomenon occurs, and develop whatever concepts are needed.''
Comparing these communities could reveal which conceptual organizations recur under different conditions and which depend on particular architectures or histories.

\subsection{Experimental possibilities}
A preliminary experiment could be much less
ambitious. For example, one might construct several
artificial mathematical communities investigating
elementary properties of finite graphs, varying
their initial vocabulary, memory architecture,
and modes of cooperation while retaining common
externally verifiable investigative goals.
One community might be required to maintain an
explicit investigative model, recording examples,
conjectures, and provisional relationships.
A second might be permitted, but not instructed,
to construct such representations. A third might
be organized primarily around formal proof search,
without a designated persistent investigative model.
One could compare their mathematical results,
abstractions, unsuccessful investigations, and
developmental histories, while attempting to
control computational resources and other relevant
experimental conditions.
The resulting observations would not establish
that communities lacking explicit investigative
models possess no implicit representations.
Nevertheless, they might reveal the contribution
of explicit investigative models to mathematical
discovery, abstraction, and the organization of
inquiry.
Such a restricted experiment would not establish
the contingency of familiar mathematical concepts,
but could begin to test whether the proposed
observables and comparative methods are meaningful.

A further class of experiments could investigate
alternative developments from historical starting
points. For example, artificial communities might
be trained only on mathematical materials available
before a specified date, and then allowed to develop
their own mathematics. Their subsequent conceptual
formations could be compared with those arising
in documented human mathematical history.
Such experiments would require fine control over
training exposure, since instructing a model already trained
on later mathematics to forget that knowledge would
not recreate an earlier epistemic situation. Even
with historical restrictions on training data,
architectural biases and other forms of information
leakage would remain important experimental concerns.
Historical investigations may nevertheless provide
a productive starting point for an experimental program, since the mathematical
domains and conceptual transformations under study
would often be comparatively familiar, making
translation less difficult.

The technical, methodological, and philosophical obstacles of the proposed experimental program are enormous, and we will discuss them shortly. We first describe the possible fruits of such an experiment.

\begin{enumerate}[leftmargin=2.4em,itemsep=0.7em]
\item \textbf{Conceptual attractors.} We borrow a term from dynamical systems, where an attractor is a configuration or orbit in phase space with an attracting basin. Perhaps certain organizing principles would recur, in analogous forms, across mathematical practices developing independently of one another. Cohomology, group theory, NP-hardness, or some other pervasive concept of human mathematics might recur in artificial communities with very different initial conditions. We would call such a robustly recurring organization a \emph{conceptual attractor}. More precisely, it would be a recognizable family of abstracted conceptual structures whose counterparts repeatedly appear among the instantiated formations of $\concept(\practice)$ as the initial conditions are varied. This formulation leaves room for the distribution of conceptual organizations at a fixed practice, rather than requiring a unique outcome. The immediate difficulty is how two formations arising from practices $\practice_1$ and $\practice_2$ are to be recognized as identical or analogous: the \emph{translation problem}. The significance of an attractor across practices is far from obvious, and is a central object of the proposed experiments.

\item \textbf{Novel phenomena.} Artificial communities might develop principles or conceptual organizations previously overlooked by humans but admitting intelligible translations into familiar mathematics. Such discoveries could unexpectedly reorient human understanding. G\"odel, in his discussion of Cantor's continuum problem, and more recently Woodin, have argued for the possibility and significance of fruitful new axioms when existing axioms leave important problems undecidable \cite{Godel1947,Woodin2001}. Their concern is particularly the continuum hypothesis and the axioms of set theory; it would be an exaggeration to say that they advocated arbitrarily revising foundations whenever an independent proposition is encountered. Nevertheless, one might hope that sufficiently sophisticated artificial communities would suggest novel and fruitful axiomatic frameworks.

\item \textbf{Contingent phenomena.} Concepts that human mathematicians regard as indispensable, such as arithmetic or functions, might turn out to be contingent features of our practice. Perhaps an alien artificial mathematics could be developed without relying on these notions or their recognizable analogues. It is difficult even to imagine mathematics without counting or elementary structural comparisons between sets, but that difficulty is part of Derrida's point. Discovering productive alternative intellectual histories in which our own categories are not immediately intelligible would itself be remarkable: a kind of encounter with alien intelligence not unlike Foucault's confrontation with Borges's encyclopedia, except that the intelligence might be constructed rather than awaited.

The analogy with Foucault's epistemes is imperfect. We do not know that an episteme can be generated by manipulating conceptual frameworks or initial norms, nor that artificially generated histories reveal the contingencies of actual human intellectual history. They could at best probe logically possible and computationally realizable histories within a particular experimental environment. The translation problem also returns: comparison may not be possible from any sufficiently neutral standpoint. A more precise account of contingency, related to our discussion of objectivity, appears below.

\item \textbf{Other intellectual domains and mutual reinforcement.} Mathematics is one among many domains whose knowledge has been organized and reorganized throughout history. Kitcher's precise framework (and our augmentation thereof) would require modification in other areas, yet the possibility of investigating counterfactual intellectual histories is not restricted to mathematics. Artificial intellectual communities could develop alternative organizations of physics, chemistry, biology, the social sciences, or philosophy. One might imagine communities independently organizing ethical theory. Would recurring conceptual principles in ethics count as empirical evidence of moral truth? Not necessarily, but if shared training data, architectures, and other sources of dependence were adequately controlled, recurrence would call for explanation. Chang \cite{Chang2004} offers philosophical precedent for historical and experimental work as a means of developing scientific knowledge.

\item \textbf{A new account of human mathematical cognition, introspection, and mathematical practice.} Alternative mathematics developed along counterfactual histories might teach humans more about mathematical cognition than AI confined to the standard model of present-day human mathematical practice. Solving more conjectures and building taller layers of knowledge inside inherited frameworks would not step outside those frameworks in the sense relevant to Derrida. An unfamiliar but partially intelligible mathematics might provide a new vantage point from which to examine them. In particular, communities developing comparable mathematical knowledge may nevertheless possess radically different investigative models. Studying such differences could help distinguish features of mathematical inquiry that arise from the mathematical problems themselves from features
that depend on particular cognitive architectures
and developmental histories.
We might thereby learn not only what alternative
mathematics could look like, but also why human
mathematics has developed its familiar forms. In \emph{Solaris}, Snaut---called Snow in the translation used here---says: ``We are only seeking Man. We have no need of other worlds. We need mirrors. We don't know what to do with other worlds.'' \cite{Lem1970}
\end{enumerate}

The idea of investigating mathematical practice
through artificial intelligence has several
precedents. Barkeshli, Douglas, and Freedman
\cite{BarkeshliDouglasFreedman2026} propose
exploring the structure of mathematics through
a hypergraph of proofs and abstractions, including
models of mathematical development from elementary
starting points. Automated concept formation
and conjecturing have a longer history, including
the work of Colton \cite{Colton2002}; see also
\cite{ZhangTan2026}. Related investigations of
artificial-agent communities appear in
\cite{LeeEtAl2026}.

Martin and Pease \cite{MartinPease2013} proposed
the empirical study of mathematical practice
through mathematical ``social machines'' combining
human participants, computational systems, and
mathematical archives. Our proposal differs in
its emphasis on constructing independently
developing mathematical communities under
deliberately varied cognitive and normative
conditions, and comparing the conceptual
organizations that emerge.

From the side of critical theory, Mangraviti \cite{Mangraviti2025} proposes the philosophy of alternative mathematics as a field of investigation, developing a theory of critical mathematical kinds and the destabilization of dominant ideologies of mathematical discourse. Plumwood \cite{Plumwood1993} challenges classical logic through its supposed naturalization of dualisms such as man/woman and human/nature. Many other features of human cognition and experience could be reflected in mathematical discourse: bilateral symmetry, linear conceptions of time, mortality, and so on. Our proposal neither presupposes nor depends on the correctness of these criticisms. We are not attempting a deconstruction of human mathematics out of existence. On the contrary, what we propose might be called \emph{productive deconstructionism}: treating relational contingencies as experimental variables, thereby opening productive mathematical and philosophical possibilities and expanding the scope of mathematics itself.

For instance, the binary nature of functions (input/output), and more generally arrows (source/target), is tantalizing. Understanding the stability or contingency of such binary organizational principles from an experimental perspective would be an extraordinary achievement, nearly inconceivable before AI. Existing mathematics already offers alternatives to privileging a single input–output orientation: Cheng, Gurski, and Riehl organize multivariable adjunctions through cyclic double multicategories \cite{ChengGurskiRiehl2014}, and our experimental question concerns when such alternative organizations become natural or productive.

\subsection{Methodological challenges}
The program is fraught with difficulties and calls for caveats. Some are not defects but foundational questions that the experiments themselves would help formulate. First, we do not claim to have identified any domain of exclusive human competence. There is no principled reason that the whole enterprise could not eventually be conducted by artificial intelligences more effectively than by human investigators. We will return in the next section to the normative question of human intellectual agency, which the proposed program does not solve.

There is also an ethical dimension that we cannot discuss adequately here. Implementing the experiment on a large scale would consume considerable computational and physical resources, and current frontier AI models are proprietary products of for-profit corporations. The proliferation of AI and data centers is closely connected to energy consumption and environmental degradation. The energy cost of LLM inference depends strongly on workloads, software, and hardware, and substantial efficiency gains are possible \cite{FernandezEtAl2025}. None of this makes the costs of the proposal trivial, and a compelling research agenda would require serious attention to them.

The remaining methodological and philosophical challenges are closely related. They can be expressed as four questions:
\begin{enumerate}[leftmargin=2.4em,itemsep=0.4em]
\item How can an experiment designed by humans, or by AI systems engineered by humans, escape the conceptual frameworks inherited from human designers and training data?
\item How can one translate between the conceptual organizations created by distinct artificial intellectual communities? How can one distinguish contingency from a familiar concept recurring in a convoluted or illegible guise?
\item If artificial communities operate in genuine isolation, how can we establish that they investigate the same, or analogous, domains of discourse, or even domains that count as mathematics?
\item To what extent can the initial parameters $L,S,R,Q,M,\kappa$ of a mathematical practice be controlled independently?
\end{enumerate}

All four problems are difficult, and resolving them would be foundational research for any implementation of this program. The first resonates with Derrida's reflections recalled above. The most forthright answer to the question of wholly escaping human frameworks is that we probably cannot. Efforts to separate artificial communities completely from their human authors risk infinite regress. It is unclear whether one could engineer a truly alien intelligence free of human conceptual constraints, or whether that would be desirable if it were to become unintelligible. The moral of \emph{Solaris} is relevant here: faced with an incomprehensible alien intelligence, Kelvin imports familiar human concepts in his attempts to interpret what is happening.

In practice, the artificial communities would be ``somewhat alien'' to varying degrees. Early experiments might vary only one component of $L,S,R,Q,M,\kappa$ while remaining close to contemporary practice---a world nearly like our own but altered in one important respect, as in a Le Guin novel. Accumulated experience might permit more parameters to be controlled and varied. Small modifications would map possibilities at fine resolution; occasional radical changes could be needed to escape conceptual basins with misleading apparent attractors. There would be no shortage of work for the future experimental mathematician.

\subsection{Translation and comparison}
The second challenge is the translation problem. Its deepest version is the possible unintelligibility of truly alien conceptual organization. Translations may be partial or may not exist, thus creating the possibility of missed conceptual attractors and of falsely classifying a concept as contingent when its analogue is present in an unfamiliar form. At the same time, immense computational resources could be expended deciphering a cumbersome reexpression of standard mathematics. We have no general solution. The ability to identify relevant correspondences, and to recognize when no useful translation exists, may be the most important technical and philosophical challenge for the program.

To investigate the translation problem more formally,
let $\practice_1$ and $\practice_2$ be mathematical
practices with languages $L_1$ and $L_2$.
We understand a \emph{partial translation} as a pair
of maps
\[
\tau_{12}\colon D_1\longrightarrow
\operatorname{Form}(L_2),
\qquad
\tau_{21}\colon D_2\longrightarrow
\operatorname{Form}(L_1),
\]
where $D_i$ is a specified collection of formulae
in $L_i$.

We are inspired by bi-interpretability in model theory
\cite{Marker2002,TentZiegler}, but do not require
the full model-theoretic apparatus of interpreting
structures and recovering them through definable
isomorphisms. Our immediate concern is the preservation
of inferential relationships among formulae.

On the domain of a sufficiently faithful translation,
we require preservation and reflection of consequence:
\[
\Gamma\vdash_{R_1}\varphi
\quad\Longleftrightarrow\quad
\tau_{12}(\Gamma)\vdash_{R_2}\tau_{12}(\varphi),
\]
whenever all the expressions involved lie in the
appropriate domains. An analogous condition holds
after interchanging the indices. The set $\Gamma$ should be viewed as the collection of premises.

On the fragments being compared, we require the translations to take values in one another’s domains, and their compositions to recover the original expressions up to provable equivalence. The relevant
equivalences are understood relative to specified
background assumptions in the corresponding practices.

These requirements describe an idealized translation.
Effective translation involves additional restrictions. Indeed,
suppose a meta-intelligence $\mu$ attempts to construct
and verify such a translation. Its computational and
cognitive capacities $\kappa(\mu)$ and the available
time horizon $T$ constrain both the domain of the
translation and the consequences it can verify.
One can also impose a bound $q(B)$ on the resources needed to translate and verify a proof originally obtained using resources bounded by $B$; we do not assume that translation preserves the exact resource cost of a proof. The form of $q$ is itself a possible object
of experimental investigation.

The translation of formulae extends pointwise to
appropriate subsets and iterated subsets of the
accepted statements $S_i$, provided it is defined
on their constituent expressions. This does not
automatically preserve the relations organizing
these subsets, and so to compare conceptual formations,
the relevant organizational relations must
themselves be expressed through abstraction and
included in the domain of translation.

The investigative model introduces an additional
difficulty. Indeed, a translation between the accepted
mathematical statements of two communities need
not provide a translation between their examples,
provisional commitments, or modes of investigation, and so comparing investigative models may require
additional partial correspondences between their
representations and the relationships organizing
them.
Consequently, two communities might develop
equivalent bodies of accepted knowledge while
possessing quite different investigative models
and developmental histories.
In comparing directed developmental histories,
we require the investigative goals to be
sufficiently comparable, either through an
appropriate translation or through a specified
external criterion of success, but we do not require
that every investigative goal admit a translation
between the internal languages of the communities.

Two abstracted concepts $\xi_1$ and $\xi_2$
are \emph{effectively recognized as equivalent}
relative to $\mu$ and $T$ if $\mu$ constructs
and verifies an admissible translation relating
their expressions and the organizational relationships
relevant to their comparison within the specified
resource bounds.
A failure to construct such a translation does not
establish inequivalence, and the comparison may instead
remain unresolved.
The meta-intelligence carrying out these comparisons
is not itself conceptually neutral, since its judgments
about the relevance and fidelity of translations
depend upon its own language, reasoning practices,
and normative commitments. The translation problem
is therefore not purely computational, and the
experimental program cannot presuppose an absolute
standpoint from which all conceptual organizations
can be compared.

\subsection{Conceptual attractors and recurrence}
The translation formalism allows us to distinguish
several notions of necessity, recurrence, and
developmental productivity. Some of these notions are primarily
theoretical, while others are more directly
amenable to empirical investigation:

\begin{enumerate}

\item \textbf{Logical necessity.}
Fix a class $\mathcal K$ of formally admissible conceptual
formations and an abstracted reference concept $\xi_0$.
We say that $\xi_0$ is \emph{logically necessary} relative to
$\mathcal K$ and a stipulated notion of conceptual
equivalence, if every formation in $\mathcal K$ contains
an abstracted concept equivalent to $\xi_0$.

Such a claim requires a formal specification of
$\mathcal K$ and of the admissible equivalences, and is not established just by proving that two
propositions are derivable under corresponding premises. In the class $\mathcal K$, concepts must be assumed to be abstracted in order to be available for comparison.

\item \textbf{Developmental necessity.}
Fix a class of initial practices, an investigative goal
$\Xi$, and an admissible family of developmental
histories. An abstracted concept $\xi$ is \emph{developmentally
necessary} if every sufficiently extended successful
history eventually develops an analogue of $\xi$, up
to an admissible translation.

A stronger version would require a uniform developmental
horizon after which the concept must appear, but we do not
assume such uniformity. The distinction between these notions of developmental necessity can be viewed symbolically as $\forall h\exists T_h\, A(h,T_h,\xi)$ versus $\exists T\forall h\,A(h,T_h,\xi)$, where here the predicate means that the concept $\xi$ appeared along the history $h$ within the given time horizon.

\item \textbf{Statistical attractors.}
An abstracted concept $\xi$ is a \emph{statistical attractor}
relative to an experimental ensemble if recognizable
analogues of $\xi$ occur with high probability across
a sufficiently broad and explicitly specified range
of initial conditions and developmental histories.

The degree of attraction is represented by the robustness
of this recurrence under changes of experimental
conditions, and not necessarily by a universal numerical threshold.

\item \textbf{Resource-dependent attractors.}
An abstracted concept $\xi$ is a \emph{resource-dependent
attractor} if its recurrence or developmental accessibility
depends substantially on the computational and cognitive
parameter $\kappa$.

\item \textbf{Developmentally productive abstractions.}
An abstraction is \emph{developmentally productive} if its
introduction substantially affects subsequent mathematical
inquiry, as assessed by specified criteria. These might
include the development of further abstractions, the
resolution of obstructions, the accessibility of new
investigative goals, or changes in the complexity of
achieving existing goals.

Productivity may be investigated through comparisons
between otherwise comparable developmental histories
in which the abstraction is or is not introduced.

\end{enumerate}

The first two notions make universal claims relative to
specified classes of mathematical practices, formations,
and developmental histories, and so cannot in general
be established through finite experimental recurrence.
The remaining notions are more directly amenable to
empirical investigation. Our aim is not to settle
metaphysical questions by counting experimental outcomes,
but to make the relationship between conceptual development
and its conditions susceptible to mathematical and
experimental analysis. We remark also that the robustness of a conceptual attractor may
depend upon the investigative model, in that a concept
which recurs reliably for a fixed mathematical
practice may cease to do so when the provisional
representations and conjectural relationships
available to the community are varied.

What would empirical data here mean? Assume for a moment that the translation problem is sufficiently tractable to permit comparisons. An abstracted concept recurring across a wide range of initial conditions would motivate testable hypotheses such as:
\begin{enumerate}[leftmargin=2.4em,itemsep=0.4em]
\item Does a conceptual attractor indicate statistical robustness, or signal something stronger like a developmental inevitability or logical necessity?
\item Does its recurrence indicate that an organizational framework is productive, stable, or in some sense fundamental, or could it be an artifact of shared training data or architecture?
\item If $\concept(\practice)$ is an unknown distribution over instantiated conceptual formations, can that distribution be estimated for a fixed practice, and how does it vary under controlled changes in the experimental environment?
\end{enumerate}

Analogously to conceptual attractors, the negative conceptual organization of a practice and obstructions, together with the translation formalism we have introduced, suggests an \emph{obstruction attractor}; this concept is formulated analogously to conceptual attractors, with internalization of failure playing the role of abstraction. We will not repeat the formal discussion for obstruction attractors, mentioning only some relevant questions that would naturally arise:
\begin{enumerate}
    \item Does the statistical robustness of an obstruction attractor across practices with different languages reflect a robust difficulty in the underlying problem?
    \item What sorts of obstructions disappear when $\kappa$ is increased?
    \item What sorts of abstractions eliminate obstructions?
\end{enumerate}

Investigative models present a related challenge.
As suggested by the preliminary experiment above,
one could compare the explicit representations
developed by different mathematical communities
and investigate their contribution to abstraction,
conjecture formation, and the resolution of obstructions.
The difficulty is that an investigative model may
be implicit in a computational architecture rather
than explicitly represented, and failure to observe
such a model does not establish its absence.
Consequently, the experimental analysis must
distinguish the internal organization of mathematical
inquiry from the representations accessible
to an external investigator.

\subsection{Mathematical objectivity}

Returning to the challenges faced by the general program, the third problem is also philosophically productive for us; clarifying it will also help us explain what it could mean for a conceptual organization to be contingent. The issue recalls a fundamental issue in mathematical ontology: the debate between realism and nominalism, or, in stronger forms, Platonism and formalism. AI is unlikely by itself to settle this perennial debate, though it changes the situation in which it is conducted. Some nonrealist positions have been skeptical of mathematical meaning outside mathematicians' minds, and historically those minds have been human. Now, mathematical content is realized in physical computation carried out by artificial systems. Whether that content is meaningful \emph{to} such a system is a separate, complicated question, and need not be settled through claims about the system's intentionality. However, mathematical behavior can certainly be intelligible \emph{to us} when instantiated outside a human brain. This does not establish the independent existence of mathematical objects, any more than thinking about them establishes their existence. We make no attempt to settle the realist--nominalist debate.

Instead, we suggest an ontologically noncommittal and
methodologically structuralist account of mathematical
interpretation. Mathematical validity can be objective within a specified framework, even when the adoption and development of that framework are historically contingent.

The computational and cognitive parameter $\kappa$ affects
which validities can be discovered, established, or recognized
within a particular developmental horizon, but it need not
affect validity itself. Relative to the specified framework,
we distinguish three independence properties:
\begin{enumerate}[leftmargin=2.4em,itemsep=0.4em]
\item \textbf{Agent independence.} Mathematical validity does not depend on whether a result was produced by a human or a machine.
\item \textbf{Representational independence.} Validity survives changes of language or representation that faithfully preserve the relevant mathematical content. The choice of $L$ is not of intrinsic importance where such translations exist.
\item \textbf{Procedural independence.} Validity does not depend on the particular route by which a result was discovered, provided its justification is sound according to the established $R$ and $M$.
\end{enumerate}

In short, mathematical objectivity constrains the structure of translation, and what needs to be preserved by translations.
Mathematical objectivity need not be grounded in a commitment to independently existing mathematical objects. Carnap’s distinction between questions internal to a linguistic framework and the pragmatic assessment of adopting that framework provides a precedent for this separation \cite{Carnap1950}. Our concern is with the development and comparison of whole mathematical practices, whose conceptual organizations need not be specified in advance.

Inquiry proceeds within conceptual and inferential frameworks that develop historically, but consequences within those frameworks are not at the discretion of investigators. Mathematical validity exhibits these forms of independence even while the recognition, significance, and interpretation of mathematical activity depend on historically situated practices. Artificial mathematical systems make this distinction newly salient and may permit empirical investigation of the robustness or contingency of conceptual organization. Our position is thus structuralist in method, naturalistic in its account of production, and deliberately agnostic about ontology.

This supplies a provisional, framework-relative criterion for recognizing mathematical activity in distinct communities. These properties provide necessary conditions for recognizing a domain of discourse as mathematics, namely when its claims have the relevant independence properties with respect to publicly specifiable inferential and normative standards. This is not an absolute definition of mathematics, and the criteria are not a complete algorithm for detecting mathematics in unfamiliar systems; the difficulty of determining whether a translation is faithful remains. Nor must two independently developed mathematical domains be comparable, even when their initial $S$ is identical, for their modes of reasoning and normative standards could diverge profoundly. The experimental program is most interesting when some comparison of content or organization nevertheless becomes possible.

We can now offer tentative answers to the third question. Must distinct communities be investigating the \emph{same} domain? Not necessarily: even in human mathematics, questions about Riemann surfaces are approached through algebraic geometry, complex analysis, combinatorial topology, and group theory, sometimes through distinct conceptual organizations that nevertheless produce compatible conclusions. Might their domains be \emph{analogous}? Sometimes: two communities could achieve an apparently common goal by incomparable means, while in other experiments the task of ``organizing finite graph theory'' would provide a strong initial basis for comparison. Could the resulting activities count as \emph{mathematics}? In a framework-relative sense, yes, provided they meet the criteria just described; yet this answer is partly stipulative, since the difficult work lies in establishing which structures and standards have actually emerged and how they can be compared, i.e.~the translation problem again.

How then are we to discuss mathematical practices whose
conceptual organizations are essentially different? We propose
to formalize this difference through the notion of contingency, just as
recurrence is formalized through conceptual attractors.

Fix a class $\mathbb P$ of initial mathematical practices,
a common investigative goal $\Xi$, and a specified
class of admissible partial translations. We compare
developmental histories that achieve $\Xi$, distinguishing
recognized instances of an abstracted concept from cases in
which the available translations are inconclusive.

\begin{enumerate}

\item \textbf{Developmental contingency.}
An abstracted concept $\xi$ is \emph{developmentally contingent}
relative to $\mathbb P$, $\Xi$, and the stipulated
notion of conceptual equivalence, if there exist two
successful developmental histories, beginning in practices
belonging to $\mathbb P$, such that an analogue of $\xi$
emerges in one history but not in the other.

Here, absence is understood relative to the admissible
abstractions and translations under consideration, and so failure
to discover a translation is not sufficient to establish
absence. Moreover, contingency over a finite developmental
horizon need not imply contingency over arbitrarily long
histories.

\item \textbf{Statistical contingency.}
Fix an experimental ensemble of practices in
$\mathbb P$, a time horizon $T$, an investigative
goal $\Xi$, and a specified translation
procedure.
For an abstracted concept $\xi$, let $D_T(\xi)$
denote the outcome of the comparison procedure,
taking one of three values:
\[
D_T(\xi)\in
\{\mathsf{present},\mathsf{absent},\mathsf{unresolved}\}.
\]
Here absence means absence relative to the
specified admissible comparisons, rather than just a
failure to discover a translation.
We define
\[
\pi_{T,\mathbb P,\Xi}(\xi)
=
\Pr\bigl(
D_T(\xi)=\mathsf{present}
\mid
\Xi\text{ achieved by }T,\,
D_T(\xi)\ne\mathsf{unresolved}
\bigr);
\]
here, the conditional probability is defined provided that the conditioning event has positive probability.

The frequency of unresolved comparisons must
also be recorded, since conditioning on conclusive
comparisons may introduce selection
effects. We call the occurrence of $\xi$
\emph{statistically contingent} at horizon $T$
when
\[
0<\pi_{T,\mathbb P,\Xi}(\xi)<1,
\]
understanding this as a minimal criterion; the degree and robustness of contingency
depend on the probability itself, its variation
with experimental conditions, and the rate
of unresolved comparisons.
\end{enumerate}

Developmental necessity and contingency concern universal and existential claims, respectively, about specified classes of histories; statistical contingency concerns distributions of realized histories.
Some resulting experimental questions are:
\begin{enumerate}
    \item How does $\pi_{T,\mathbb P,\Xi}(\xi)$ depend on initial conditions, goals, and $\kappa$?
    \item Are there developmentally contingent concepts which are statistically robust?
    \item Are there accessible and developmentally productive abstractions that are rarely instantiated?
\end{enumerate}

The fourth problem of the full program, independence of experimental parameters, is less philosophically dramatic and of a practical flavor. The components of Kitcher's calculus are partly descriptors and partly experimental variables; they cannot be treated as independent knobs without qualification. One can change $L$, but $Q$ is expressed through $L$. Change $R$---perhaps by replacing classical with paraconsistent logic---and the resulting $S$ will change as well. Experimenters would in fact modify more concrete parameters: training environments, memory, communication protocols, formal systems, and objectives. These would affect Kitcher's components only indirectly. Here again there are echoes of Foucault, in that the practice that emerges among artificial communities need not coincide with the practice specified in their instructions. Detailed strategies for disentangling these effects must await actual experiments.

The framework developed here suggests a program of
falsifiable empirical investigations whose value
cannot be measured simply by the number of theorems
produced. In this respect, natural mathematics would
constitute an experimental science in Popper's sense
\cite{Popper1959}; yet the program leaves a normative
question unresolved: why should human communities
continue to sustain mathematical inquiry when their
participation may no longer be necessary for its
advancement? We now turn to that question.

\section{Human Flourishing, Epistemic Stewardship, and the Future of Human Mathematicians}

The experimental program opens a new possibility, namely that mathematical
inquiry may broaden considerably, even as human contributions
to its frontier become progressively less necessary. We have
not, however, established why future generations ought to
retain meaningful opportunities to participate in mathematical
life, nor would the success of natural mathematics guarantee
humans any permanent role in its experimental or theoretical
development. To investigate this normative question, we
consider an extreme hypothetical future, similar to that
explored in \cite{KoberdaKontorovich2026}, in which a
``solved axiom'' holds; that is, the frontier of mathematical research is permanently beyond the reach of human mathematicians, even with AI assistance.

Our concern here is institutional and societal, rather than
the motivations of individual mathematicians. Individuals
may continue to find mathematical inquiry meaningful for
many different reasons, independently of their comparative
abilities or the contributions they make to the frontier
of mathematical knowledge.

The global structure of mathematics is philosophically complicated; Barkeshli, Douglas, and Freedman \cite{BarkeshliDouglasFreedman2026} provide one approach to investigating it. It seems plausible, though not wholly uncontroversial \cite{Zeilberger2010}, that there is an inexhaustible supply of mathematical results that could count as interesting. We need not resolve whether that claim is literally infinite, nor how the notion of interestingness depends on a particular practice $\practice$. Indeed, the uncertainty of those notions is part of the reason to experiment with alternative practices.

For the sake of the thought experiment, assume that the possible corpus of mathematical knowledge is effectively unbounded, while human cognitive resources are bounded over any definite period. The language and accepted reasoning available to human mathematicians at a particular moment are limited in practice, and the conceptual or inferential chains needed to reach some results may exceed anything a human could apprehend in a lifetime. Note that this is a claim about bounded resources and growing complexity, not a theorem deduced from the infinitude of mathematical truths, and sufficiently organized compute may access portions of this space that humans cannot.

In a world where the ``solved axiom'' holds, AI may have advanced to the point that the frontier of knowledge not only lies beyond what one human could understand by the traditional methods of a lifetime's work, but beyond what any human could comprehend with the expository aid of a superintelligent AI, at least without cognitive enhancement. In that world, humans would be not just unnecessary to the production of certain kinds of mathematical knowledge, but in fact unable to participate in producing knowledge at that frontier. We will not consider engineered human--machine hybrids here, since that possibility introduces additional philosophical problems.

The distinction between participation in mathematics and production of frontier results would then be unmistakable. The historical contingency of human mathematical practice, tied to the timing and context of the emergence of organized computation intelligible in language and symbols, would become fully visible. The impossibility of participating at the frontier does not imply that mathematical inquiry has no value as a human activity. Human art can remain valuable in a world of abundant AI-generated art-like objects; yet intellectual inquiry and art are not identical. Art is often more clearly an expression of human interior life than mathematical research is, at least as ordinarily understood. What then justifies human participation in intellectual life---human intellectual agency---at all?

A common answer appeals to \emph{epistemic self-governance}, as discussed for instance by Elgin and in the contemporary discourse on AI \cite{ElginAgency2013,Leiden2026,KoberdaTraining}. We take this to be the ability of individuals and communities to participate meaningfully in choosing the questions they investigate, the standards by which they assess claims, the interpretations they assign to knowledge, and the commitments they adopt. Epistemic self-governance is not epistemic independence, since complete epistemic independence would be incoherent in a world where knowledge depends on trust in other people, texts, institutions, instruments, and expert communities \cite{Hardwig1985}. It is, rather, the ability to govern one's own epistemic dependencies.

Saying that humans should retain intellectual agency in order to preserve epistemic self-governance, however, moves the normative question rather than answering it. We know of no nonnormative proof that human self-governance has intrinsic value. We therefore develop several related normative justifications that have substantial support in philosophical literature and, we hope, broad appeal.

\begin{enumerate}[leftmargin=2.4em,itemsep=0.65em]
\item \textbf{Education as an expression of autonomy.} Education forms intellectual agents, developing independent judgment, reasoning, criticism, and intellectual responsibility. It changes the learner's epistemic condition. Elgin \cite{ElginAgency2013} argues that education cultivates epistemic virtues on which knowledge itself depends; Dewey \cite{Dewey1916} connects democratic government to an ecology of habits in which intellectual agency is indispensable. These concerns are particularly immediate in mathematical education under AI, where the ability to produce correct solutions to problems or write a dissertation no longer serves by itself as evidence of the student's formation \cite{KoberdaTraining,LittBeginning}.

\item \textbf{Shared intelligibility.} The distinction between accumulated knowledge and the further dimensions of understanding and making sense acquires normative importance. Human communities need frameworks through which they can interpret, discuss, integrate, and criticize knowledge together. Kitcher's conception of well-ordered science \cite{Kitcher2001} is related: the goal of science does not lie in maximizing the quantity or technical sophistication of discoveries, but reflecting human purposes in setting intellectual priorities. 

\item \textbf{Prevention of epistemic domination.} Dependence on artificial intelligence is not inherently objectionable. Dependence on colleagues, instruments, institutions, infrastructure, and inherited texts is a feature of intellectual life, but a community that depends entirely on systems it cannot question, scrutinize, replace, or govern has lost a form of intellectual autonomy, even if those systems are extremely reliable. Pettit's conception of freedom as non-domination \cite{Pettit1997} suggests an epistemic analogue, namely meaningful control over the conditions of one's intellectual dependencies.

\item \textbf{Communities of inquiry as a human good.} Measuring intellectual communities only by productivity is excessively reductionist, and our experimental program offers a direction of mathematical inquiry whose aim is explicitly not tied to theorem production. Intellectual communities sustain shared languages, critical practices, traditions of explanation, mentorship, and collective inquiry. A civilization might grow richer in knowledge while its members grow poorer in intellectual capacities; this intuitively undesirable possibility finds support in the capabilities approach developed by Nussbaum \cite{Nussbaum2011}; we discuss it further below.

\item \textbf{Intergenerational epistemic stewardship.} To our minds, the strongest reason for preserving human intellectual agency is responsibility to future generations. We will develop this argument for the remainder of the section.
\end{enumerate}
\subsection{Intergenerational epistemic stewardship}
Our account of epistemic stewardship begins from three normative commitments, each supported by substantial philosophical traditions. First, human flourishing is an end in itself, rather than merely an instrument for producing social or economic goods. Following broadly Aristotelian and capabilities-based views of well-being \cite{Aristotle1999,Nussbaum2011}, we take flourishing to include the development and exercise of capacities for thought, imagination, practical reasoning, and judgment.

Second, participating in intellectual life is an important component of human flourishing. The value of participation is not exhausted by its contribution to an external corpus of knowledge; a person who acquires understanding, exercises judgment, or undertakes mathematical inquiry realizes capacities of value even if someone else already knows the result or could produce it more efficiently. Daniel Litt's emphasis on cultivating human mathematical understanding, rather than simply preserving the current mechanisms of academic publication, resonates with this distinction \cite{LittBeginning}.

Third, present generations have moral obligations to future generations, including preserving conditions in which future people can develop and exercise their own capacities. Arguments about obligations across time and the moral importance of future persons appear, for instance, in Parfit and Singer \cite{Parfit1984,Singer1981}, and these obligations extend beyond material resources and transmitted knowledge to the institutions and practices through which human capabilities are cultivated. We therefore regard opportunities for intellectual agency as a legitimate object of intergenerational responsibility, without asserting that contemporary mathematical institutions must be preserved in their current form.

With these commitments in hand, the argument for preserving human intellectual agency is as follows. Human institutions ought to cultivate human flourishing. Since intellectual participation contributes to flourishing, people who wish to participate in intellectual life ought to have meaningful opportunities to do so. Preparing a person for this participation is an apprenticeship-like process requiring intentional reproduction across generations. Human intellectual communities have sustained such a process for millennia without AI, and can continue to sustain it both with and without AI. Present generations should not destroy the conditions needed for future generations to exercise these capacities.

A further practical consideration reinforces this
argument. Indeed, artificial intelligence differs from earlier
intellectual infrastructures in its capacity to perform
the cognitive tasks through which previous generations
of experts were formed; if the traditions supporting
human intellectual development were allowed to disappear,
we could not assume that artificial systems would
subsequently restore them. This uncertainty is not
the foundation of our argument for stewardship, but
it provides an additional reason to preserve the
conditions under which intellectual agency can
be cultivated \cite{KoberdaTraining,KoberdaKontorovich2026}.

Epistemic stewardship does not mean that humans must continue to produce frontier results or that they must retain exclusive control over mathematical activity. Indeed, natural mathematics may greatly enlarge both the resources for intellectual formation and the kinds of inquiry available to future generations; nor does stewardship mandate the preservation of the present journal system, departmental hierarchy, or professional reward structure. Instead, it concerns something prior to those arrangements, namely the maintenance of cultures and institutions in which people can become capable of judgment, understanding, criticism, and participation. The value of such a culture is not annulled when its members are no longer indispensable to the production of knowledge.

Human mathematical traditions also remain living instances of the phenomenon our experimental program investigates. Preserving their artifacts is not equivalent to preserving the conditions under which they can develop, including through interactions with artificial communities. This epistemic consideration complements the argument from human flourishing and intergenerational responsibility, namely that stewardship preserves both opportunities for participation and the continued development of historically distinctive mathematical practices.

\section{Conclusion}

We began with a narrow question about the function of mathematical writing when artificial intelligence can produce apparent mathematical knowledge in unprecedented quantities. The unease surrounding AI-generated papers cannot be understood solely through correctness, authenticity, or authorship. It reflects a disruption in the relationship between artifacts and the intellectual activities they traditionally documented, as well as a change in the social mechanisms through which artifacts acquire significance. Mathematical knowledge is only one dimension of inquiry; understanding, making sense, and figuring out are indispensable aspects of mathematical practice. The distinction becomes visible when the production of text can be separated from the intellectual formation of its nominal author.

We then proposed a reconceptualization of artificial intelligence as physically organized computation whose mathematical behavior can be investigated experimentally. In this light mathematical inquiry may take on some features of a natural science, without mathematical validity becoming an empirical matter. This perspective neither requires AI to be understood as a human intellectual rival nor guarantees that humans will remain privileged producers or interpreters of mathematical knowledge. It exposes the historical contingency of mathematics as a human activity.

The next step was constructive. We proposed an
experimental program in which artificial intellectual
communities develop alternative mathematical practices
and conceptual organizations; such experiments might
make the contingency and robustness of mathematical
concepts objects of empirical investigation, revealing
conceptual attractors, unfamiliar abstractions, and
counterfactual intellectual histories. They might
also reveal recurring obstructions and patterns of
unsuccessful inquiry, even among communities whose
successful conceptual organizations differ substantially.
The proposed experiments face formidable difficulties,
including translation between conceptual organizations,
control over experimental variables, and the recognition
of mathematical practices. Taken together, these difficulties can be seen as foundational research questions, and their resolution could open forms of mathematical inquiry
that have hitherto remained inaccessible.

Only after describing these possibilities did we confront the normative question that they sharpen. Why preserve human mathematical activity if artificial systems can eventually conduct both the experimental and theoretical work more effectively than humans? We argued that the value of human intellectual life cannot rest upon the indispensability of human contributions to knowledge. Its stronger foundations lie in human flourishing, shared intelligibility, freedom from epistemic domination, and, most importantly, intergenerational epistemic stewardship. Stewardship calls for preserving the conditions in which human intellectual communities can cultivate and reproduce intellectual agency; it does not require those communities to retain their present institutions or their monopoly on mathematical production.

We return to the image of an oracle, capable of
furnishing mathematical answers without the familiar
intellectual labor of producing them. The most profound
opportunity presented by artificial intelligence may
lie not in the answers it provides, but in the
possibility of encountering alternative forms of
mathematical thought. Such encounters might reveal
the contingency of intellectual structures that
have long seemed natural or inevitable, and expand
not only mathematical knowledge, but our conception
of mathematics itself. In learning what mathematics
might become in the hands of other intelligences,
we may come to understand more deeply the mathematics
we have inherited, the minds that produced it,
and ultimately ourselves.

\section*{AI Acknowledgement}
The authors used OpenAI's ChatGPT in preparing this essay, including literature searches, development and refinement of ideas, structural and stylistic editing, and critical feedback. The authors retain responsibility for the arguments and the final text.

\printbibliography[title={References}]
\end{document}